\documentclass[11pt]{pspum-l}

\scrollmode
\usepackage{amsmath}
\usepackage{amsthm}
\usepackage{amssymb}
\usepackage{amsfonts}
\usepackage{epsfig}
\usepackage[matrix,arrow,curve]{xy}

\newtheorem{theorem}{Theorem}[section]
\newtheorem{pro}[theorem]{Proposition}

\newtheorem*{teo}{Theorem}

\theoremstyle{definition}
\newtheorem{definition}[theorem]{Definition}

\theoremstyle{remark}
\newtheorem{remark}{Remark}

\newcommand{\cotanh}{\mathop{\mathrm{cotanh}}}
\newcommand{\<}{\langle}
\renewcommand{\>}{\rangle}

\renewcommand{\Im}{\operatorname{Im}}

\newcommand{\nc}{\newcommand}
\nc{\Symm}{{\on{Sym}}}
 \newcommand{\q}{{\quad}}

\newcommand{\diag}{\mathrm{diag}\,}

\newcommand{\beq}{\begin{equation}}
\newcommand{\eeq}{\end{equation}}
\usepackage{hyperref}
\usepackage{cleveref}

\begin{document}
\author{Michael Semenov-Tian-Shansky}
\title[Quantum Toda Lattice:  a Challenge for Representation Theory]{Quantum Toda Lattice: \\a Challenge for Representation Theory}
\address{Institut de Math\'ematiques  de Bourgogne, UMR 5584 CNRS,
                                                             Universit\'e de Bourgogne
                                                             B.P. 47870
                                                             21078 Dijon Cedex France,
                                                             and Steklov Mathematical Institute, St.Petersburg, Russia.}

\email{semenov@u-bourgogne.fr}
\thanks{This work was supported in part by the Russian Foundation for Basic Research (Grant No. 18-01-00271).}
\date{}

\begin{abstract}
Quantum Toda lattice may solved by means of the Representation Theory of semisimple Lie groups, or alternatively by using the technique of the Quantum Inverse Scattering Method. A comparison of the two approaches, which is the purpose of the present review article,  sheds a new light on Representation Theory and leads to   a number of challenging questions.
\end{abstract}

\maketitle

\section{Introduction}
In May 1978 L.D.Faddeev delivered a program talk in his research seminar at the Steklov Institute in Leningrad. His starting point was the breaking news he had just brought from the US: this was the stunning discovery by B.Kostant of an unexpected link between Integrable systems and the famous Orbits Method in Representation theory of Lie groups \cite{K}. Kostant applied his approach to one of the well-known examples of integrable systems, the Toda lattice. In classical case, this approach yields a natural Lax pair for the Toda lattice together with the Poisson structure and reduces its solution to a standard factorization problem in the associated semisimple Lie group. In quantum case, it provides commuting quantum integrals of motion and gives explicit formulae for the Toda wave functions in terms of specific matrix coefficients of irreducible representations of the same group. One key message of Kostant's approach was that classical integrable systems admit natural quantum analogs which are still integrable and exactly solvable. The big challenge, which was the main content of Faddeev's talk, was to create an effective quantization technique which would allow to bring into play other examples of quantum integrable systems, in particular, infinite dimensional. The alternative method   outlined in his talk was a combination of the ideas of the classical inverse scattering with the algebraic methods developed by Baxter. Within a year, this method brought about a decisive breakthrough, providing a unified and highly effective approach to a large variety of integrable models. A comprehensive story of this discovery   is described by Faddeev in \cite{QISM}. The key new topics which had been developed in the course of that remarkable year were the quantum $R$-matrices, quantum Lax operators and the algebraic Bethe ansatz. It gradually became clear that the algebra involved in this new method (which is now known under the name of Quantum Inverse Scattering Method, or QISM) goes well beyond the conventional theory of Lie algebras and Lie groups.

The initial idea of Kostant which had played  the r\^ole of a catalyst in the development of QISM was somewhat eclipsed by these discoveries. While its full potential has soon become clear in application to classical integrable systems (see,  for instance \cite{RS}), its use in the quantum case remained rather limited. In the classical case, the proper generalization of Kostant's approach immediately brings into   play infinite dimensional Lie groups and Lie algebras, mainly, loop algebras  and their central extensions. The quantum counterpart of this approach would require a very precise information on their unitary  representations, which was at the time (and in fact, still is) completely lacking. The case of open Toda lattice originally studied by Kostant remained for quite a while a brilliant but somewhat isolated example\footnote{The seminal lecture of Kostant \cite{K} remained his only publication on the quantum Toda lattice. A comprehensive treatment of the spectral problem for the  Toda lattice was given by the author \cite{STS}.}.  Another example which has much in common with the open Toda lattice, is the theory of spherical functions developed by Harish-Chandra \cite{HC} as a part of his comprehensive work on representation theory  of semisimple Lie groups. As noticed by Olshanetsky and Perelomov (see \cite{OPer}) already in the early 1970's,
 spherical functions are related to another important class of quantum integrable systems, the systems of Calogero--Moser type.  Both the Toda lattice and the Calogero--Moser system share one key property, the factorized formulae for the scattering matrices which had been discovered by Gindikin and Karpelevich \cite{GinK} in the aftermath of Harish-Chandra's work on spherical functions.

 For several decades the ideas of Harish-Chandra dominated the development of the representation theory.
 As we understand it now, a new breakthrough in this classical area can be achieved by a  combination of  the old ideas of Representation Theory with those of  QISM. The application of the highly powerful method of  quantum separation of variables is of particular importance in this respect. This method, originally introduced in the pioneering works of Sklyanin \cite{Skl1}, \cite{Skl2} as a generalization of the algebraic   Bethe ansatz technique,  leads to a very strong amplification of the old factorization formulae of Gindikin--Karpelevich type. While some of the key new formulae which go in this direction have been obtained more than 15 years ago, mainly, by the former ITEP group  (D.Lebedev, S.Kharchev, A.Gerasimov, see \cite{Leb1}, \cite{Leb2}, \cite{Leb3}), they are still not widely known by experts in Representation Theory. It is important therefore to bring more attention to these results. This is the purpose of the present review.  The author believes that the understanding of these new results in their proper generality is a highly important task which will eventually reshape Representation Theory and certainly requires a lot of work. In this review we leave aside another major development also connected with the pioneering work of the ITEP group: this is the connection of the Quantum Toda lattice with Quantum Cohomology;  its exploration had started with the now famous paper  of    Givental and Kim \cite{KG} and eventually lead to a plethora of new integral formulae for the wave functions of the Quantum Toda lattice \cite{GLO}. Putting these results into the general context of  Representation Theory is another challenge which must be addressed separately.

It would have been very nice to prepare this article for the 70$^{\text{th}}$ anniversary of Boris Dubrovin, my highly esteemed colleague  and my friend for more than 40 years. It is terribly sad that instead  it will appear as a tribute to his memory.

\section{Quantum Toda Lattice and Representation Theory}\label{rep}

In order to explain the new features which are introduced by the application of the QISM techniques let us start with a brief review of Representation Theory as it has emerged in the fundamental work of Gelfand and Harish-Chandra, among many others. Let $G$ be a semisimple group (or reductive group, in order not to exclude the case $G=GL(n)$), $\mathfrak{g}$ its Lie algebra, $K\subset G$ its maximal compact subgroup.The key idea which links Representation Theory of $G$ to integrable systems is the use of Laplace operators on $G$ or on $G/K$  as a source of quantum integrals of motion  (the use of symmetric space $G/K$ is technically more simple,  as we can avoid dealing with non-elliptic operators).  The ring $\mathcal{D}(G)^G$  of Laplace operators, i.e.  $G$-invariant differential operators on $G$ is isomorphic to the center $\mathcal{Z}$ of the universal enveloping algebra $U(\mathfrak{g})$; for split real groups $\mathcal{Z}$ is also isomorphic to the ring  $\mathcal{D}(G/K)^G$ of Laplace operators on $G/K$.   Quantum integrals of motion arise as   \emph{radial parts} of Laplace operators associated with an appropriate separation of variables on $G$ or on $G/K$.  Let    $A\subset G$ be a   split Cartan subgroup, $\mathfrak{a}$ its Lie algebra. We fix an order in the system of roots of $(\mathfrak{g},\mathfrak{a})$; let $\mathfrak{u}\subset\mathfrak{ g}$ be the maximal nilpotent subalgebra generated by root spaces which correspond to positive roots.
Let $U\subset G$ be  the corresponding maximal unipotent subgroup. The group $G$ admits two standard global  decompositions, the Cartan decomposition  $G=KAK$,   and the Iwasawa decomposition $G=UAK$. Historically, the first non-trivial example of this kind  is associated with the Cartan decomposition. Let $\mathcal{C}_0=C(K\backslash G/K)$ be the space of smooth   functions on $G$ which satisfy the functional equation
 $$\phi(k_1gk_2)=\phi(g)\q  \text{for all}\q   k_1,k_2\in K.$$
 Such functions are completely determined by their restriction to $A$. Laplace operators on $G/K$ leave  $\mathcal{C}_0$ invariant and induce differential operators acting on $A$; their common eigenfunctions are called spherical functions on $G/K$. For $G=SL(N, \mathbb{R})$ we may identify $\mathfrak{a}$, the Lie algebra of $A$,  with the space of traceless diagonal matrices $x=\diag(x_1,\dots, x_N)$. Using the exponential parametrization of $A=\exp \mathfrak{a}$, the radial part of the Laplace-Beltrami operator on $G/K$ is given by
\begin{equation}\label{rad}
  -\Delta_A\phi =-\sum_k\frac{\partial^2 \phi }{\partial x_k^2}+\sum_{k<j}\cotanh(x_k-x_j)(\partial_{x_ki}-\partial_{x_j})\phi .
\end{equation}
As first noticed by Olshanetsky and Perelomov (see \cite{OPer}), by  a simple change of variables
 \begin{equation}\label{ch}\phi=\prod_{k<j}\sinh^{1/2}(x_k-x_j)\cdot \psi
 \end{equation}
we get from $\Delta_A$ the Sutherland Hamiltonian (also known as the trigonometric Calogero--Moser Hamiltonian)
\begin{equation}\label{Suth}
  H=-\Delta+\sum_{k<k}\frac{1}{\sinh^2(x_k-x_j)}.
\end{equation}
Thus the spectral theory of the Sutherland system is reduced to the theory of spherical functions on $SL(N)/SO(N)$.

The standard way to construct the eigenfunctions of Laplace operators is to express them as appropriate matrix coefficients of irreducible representations of $G$. The relevant representations are  \emph{class one representations of the principal series}. Let $M$ be the centralizer of $A$ in $K$, $B=MAU$ the minimal parabolic subgroup of $G$. We denote by $M'$ the normalizer of $A$ in $K$; the quotient $W=M'/M$ is the restricted Weyl group acting on $\mathfrak{a}$. By definition, the principal series representations are induced by 1-dimensional  representations of $B$, $\pi_\mu= ind\,_B^G \chi_\mu$, where
$$
 \chi_\mu(man)=e^{\mu(\log a)}, \;\mu\in \mathfrak{a}^*_\mathbb{C}.
$$
These representations are unitary if $\mu=i\lambda+\rho$, where $2\rho$ is the sum of all positive roots of $\mathfrak{a}$ and $i\lambda\in \mathfrak{a}^*_\mathbb{C}$ is purely imaginary.   By a slight abuse of notation we shall write $\pi_\lambda$ instead of $\pi_\mu$, etc.  By definition,  the corresponding representation space $\mathfrak{H}_\lambda$ consists of smooth functions on $G$ satisfying the functional equation
$$
\phi(bx)=\chi_{i\lambda+\rho}(b)\phi(x), \; b\in B.
$$
The group $G$ acts on $\mathfrak{H}_\lambda$ by  right translations, $\pi_\lambda (g)\phi(x)=\phi(xg)$
For $x\in G$ let $x=nak$ be its Iwasawa decomposition; we write $k=\kappa(x)$, $a=e^{H(x)}$. Put $\varphi_\lambda(x)=e^{\langle i\lambda+\rho, H(x) \rangle}$. Clearly, $ \varphi\in \mathfrak{H}_\lambda$ and  $\pi_\lambda (k)\varphi_\lambda=\varphi_\lambda$ for all $k\in K$. It is easy to see that $\varphi_\lambda$ is the only $K$-invariant function in $\mathfrak{H}_\lambda$.  Restricting functions from $\mathfrak{H}_\lambda$ to $K\subset G$ we get a model of the principle series representation in $L_2(K)$ or in $L_2(K/M)$\footnote{In an obvious way the action \eqref{prin} commutes with the action of $M$ via left translation and hence projects down to a well-defined action in $L_2(K/M)$. This latter action is already irreducible. };  the action of $G$ in $L_2(K)$ is given by
\beq\label{prin}
\pi_\lambda (g)\phi(k)=e^{\<i\lambda+\rho,H(kg)\>}\phi(\kappa(kg)).
\eeq
In this model the subgroup $K$ acts simply by right translations and the $K$-invariant vector is a constant function   $\mathbf{1}\in L_2(K)$.

All principal series representations are infinitesimally irreducible, that is, the restriction of Laplace operators $z\in \mathcal{Z}$ to $\mathfrak{H}_\lambda$ is a multiplication operator which depends polynomially on $\lambda$; in this way we get a homomorphism
\beq \label{HC_hom} \gamma: \mathcal{D}(G/K)^G\to P(\mathfrak{a}^*).
 \eeq
 By a famous Harish-Chandra theorem  $\gamma$ is actually an isomorphism onto its image which coincides with the algebra of $W$-invariant polynomials. The matrix coefficient of the principal series representation between two spherical vectors
 \begin{equation}\label{spher}
 \phi_\lambda(g)=\langle \pi_\lambda (g) \varphi _\lambda, \varphi _\lambda \rangle_{\mathfrak{H}_\lambda}=\int_{K/M}e^{\langle i\lambda+\rho, H(gk) \rangle}\,dk
\end{equation}
 is a $K$-biinvariant function on $G$. By construction, $\phi_\lambda$ is a common eigenfunction of the ring $\mathcal{D}(G/K)^G$. This is the famous Harish-Chandra formula for the spherical function on $G$;  the change of variables \eqref{ch} transforms it into the eigenfunction of the Sutherland system.

 In order to get the open Toda Hamiltonian
 \beq\label{T_open}
  H=-\Delta+\sum_{k=1}^{N-1}e^{x_{k+1}-x_{k}}, \q x\in \mathbb{R}^N, \q \sum_k x_k=0.
  \eeq
  we need to use, instead of the spherical vector, the \emph{Whittaker vector} in the representation space (the discovery of the key r\^ole of the Whittaker vectors in this context is due to Kostant \cite{K}). Let $\mathfrak{v} \subset \mathfrak{g}$ be the maximal nilpotent subalgebra spanned by all root space vectors which correspond to negative roots of $(\mathfrak{g}, \mathfrak{a})$; it admits a natural $\mathbb{N}$-grading (the \emph{principal grading}) in which all elements of degree 1 correspond to simple roots;  its commutant $[\mathfrak{v},\mathfrak{v}] $ is generated by elements of degree $\geq2$ and hence by root space vectors which correspond to non-simple roots. The characters (1-dimensional representations) of $\mathfrak{v}$ identically vanish on  $[\mathfrak{v} ,\mathfrak{v} ] $ and hence are completely specified by their values $c_\alpha= f(e_{-\alpha})$ on root space vectors which correspond to simple  roots. A character $f$ is called non-degenerate if $c_\alpha\neq0$ for all $\alpha\in P$.
   We define a unitary character $\chi_f$ of the corresponding maximal unipotent subgroup $V$ by
$$
\chi_f(e^X)=e^{if(X)}.
$$
Let $\mathcal{H}$ be a representation space of $G$, $T$ the corresponding representation.
A vector $w_f\in \mathcal{H}$   is called a \emph{Whittaker vector} if
$$
T(v)w_f=\chi_f(v)\cdot w_f \q\text{for all}\q v\in V.
$$
A natural model of the principal series representations adapted for the study of Whittaker vectors is associated with the Bruhat decomposition, which defines an embedding  of $V$  onto an open dense cell in $G/B$. In this model   $\mathfrak{H}_\lambda$ is identified with $L_2(V)$, the subgroup $V\subset G$ acts by right translations  and the Whittaker vector (which belongs to a suitable completion of $L_2(V)$) is a simple exponential, $w_\psi(v)=e^{if(\log v)}$. Equivalently,  we can define an element $w_{\psi, \lambda}\in\mathfrak{H}_\lambda$ by setting
$$
w_{\psi, \lambda}(bv)=\chi_\lambda(b)\psi(v),\; b\in B, \;v\in V.
$$
This formula defines $w_{\psi, \lambda}$ on the open dense subset $BV\subset G$ and hence by continuity on the entire group $G$.

The\emph{Whittaker function} is the matrix coefficient
\begin{equation}\label{wh}
  W_\lambda(g)=\<T_\lambda(g)w_f,\, \mathbf{1}\>_{\mathfrak{H}_\lambda}.
\end{equation}
Definition \eqref{wh} immediately implies a formal  integral representation
\begin{equation}\label{whit}
W_\lambda(g)=\int_Ve^{\<i\lambda+\rho, H(vg)\>}{\psi(v)}\,dv
\end{equation}
which is similar to \eqref{spher}; regularization and analytic continuation of \eqref{whit} can be made in the same way as for  \eqref{spher}.

The Whittaker function  $W_\lambda$ satisfies functional equation
\begin{equation}\label{whf}
W(kgv)=\chi_f(v)W(g) \q \text{for all }\; k\in K, v\in V;
\end{equation}
hence such functions may be regarded as functions on $X=G/K$ and are completely determined by their restriction to $A\cdot x_0\simeq A=\exp \mathfrak{a}$. The radial part of the Laplace--Beltrami operator for Whittaker functions in exponential parametrization is precisely the generalized open Toda Hamiltonian\footnote{The standard open Toda Hamiltonian \eqref{T_open} corresponds to the root system $A_{N-1}$ with all coefficients  $ c_\alpha=1$. }
\begin{equation}\label{Tod_G}
  H=-\Delta+\sum_{\alpha\in P} c_\alpha^2 e^{-2\alpha(x)}, \q x\in \mathfrak{a}.
\end{equation}
 A character $f$ is called \emph{non-degenerate} if all coefficients $c_\alpha=f(e_\alpha)$ are non-vanishing. After an appropriate change of variables on $\mathfrak{a}$, $x\mapsto x+x_0$, we can assume that all $c_\alpha=1$. Degenerate characters give rise to generalized Toda lattices with some of potentials omitted; such Hamiltonians are actually associated with root systems obtained by removing some vertices from the original Dynkin diagram.
 For $\mathfrak{g}=\mathfrak{gl}(n)$ we get the original Toda Hamiltonian \eqref{T_open}. By construction, Whittaker functions are its eigenfunctions (in fact, they are common eigenfunctions of the whole ring $D^G(X)$ of invariant differential operators on $X$). Whittaker functions decay very rapidly along any ray outside the positive Weyl chamber in $\mathfrak{a}$; in the positive Weyl chamber their asymptotics is given by a linear combination of plane waves with coefficients which describe scattering of wave packets on its walls.

 \begin{remark} An alternative definition of Whittaker functions makes use of a \emph{pair} of Whittaker vectors associated with two characters $\chi_f, \chi_{f'}$ of two opposite unipotent groups $U, V$; functional equation \eqref{whf} is replaced by
 \begin{equation}\label{wh2f}
W(nav)=\chi_f(n)\chi_{f'}(v)W(a), \q n\in U, v\in V.
\end{equation}
This definition leads to the Hamiltonians
\begin{equation}\label{Tod_G-ind}
  H=-\Delta_{a}+\sum_{\alpha\in P} c_\alpha c'_\alpha  e^{-\alpha(x)}, \q x\in \mathfrak{a},
\end{equation}
which depend on two sets of coefficients $c_\alpha=f(e_\alpha),\, c'_\alpha=f'(e_{-\alpha})$. The potentials in these Hamiltonians are not necessarily positive, which leads to some subtleties in the corresponding spectral theory (the Hamiltonian \eqref{Tod_G-ind} is not essentially self-adjoint on the space of rapidly decreasing functions on $\mathfrak{a}$ because of the possible ``rapid escape to infinity''). If we choose $c_\alpha=f(e_\alpha),\, c'_\alpha=f'(e_{-\alpha})=1$ for all $\alpha\in P$, we get the same Whittaker functions as in \eqref{wh} (up to a simple rescaling $x\mapsto 2x$).
\end{remark}

Following the general ``Harish-Chandra philosophy'', the spectral theory of Hamiltonians \eqref{Tod_G-ind}, \eqref{Suth}, or, more generally,  of the whole ring $\mathcal{D}^G(G/K)$ of invariant differential operators, is based on the study of ``wave packets''
\beq\label{wpack}
\Phi_{a}(g)=\int_{\mathfrak{a}^*}\<\pi_\lambda(g)a(\lambda), \mathbf{1}\>\, d\lambda,
\eeq
where $\mathbf{1}\in L_2(K/M)$ is the spherical vector and the amplitude $a(\lambda)\in L_2(\mathfrak{a}^*\times K/M)$. Clearly, $\Phi_{a}$ depends only on the coset class $x=gK$ and hence is a function on $X=G/K$. To prove that the wave packets lie in $L_2(X)$ one has to study the asymptotics of the matrix coefficients $\<\pi_\lambda(x)a(\lambda), \mathbf{1}\>$
 A key part of this study is the r\^ole of the Weyl group in its double guise of the geometric symmetry group of the root system and of the Galois group associated with the algebraic extension $P(\mathfrak{a}^*)^W\subset P(\mathfrak{a}^*)$ which underlies the definition of the Harish-Chandra isomorphism \eqref{HC_hom}. The symmetry with respect to $W$ implies that each point $\lambda\in\mathfrak{a}^*$ of the spectrum comes along with all its images $\lambda^w\in\mathfrak{a}^*$,   $w\in W$, and, moreover, all representations $\pi_{\lambda^w}$,  $w\in W$, are equivalent. This equivalence is explicitly realized by intertwining    operators introduced by Gelfand under the name of operators of horospheric automorphisms; they are closely related to the scattering theory for the wave packets \eqref{wpack}\footnote{For more details on the r\^ole of scattering theory in connection with the representation of the principal series and the intertwining    operators see \cite{STS1}.}.

 To define the  intertwining    operators we need some notations. Let $\Delta_+, \Delta_-,$ be the sets of positive and negative roots of $(\mathfrak{g, a})$. For $s\in W$ put $\Delta(s)=  \Delta_+  \cap s\Delta_-$. Put
 $$
 \mathfrak{v}_s=\bigoplus_{\alpha\in\Delta(s)}\mathfrak{g}_{-\alpha}, \; V_s=\exp\mathfrak{v}_s.
 $$
 Let $dv_s$ be the Haar measure on $V_s$. For $s\in W=M'/M$ we shall choose a representative in $M'$ which will be denoted be the same letter.  Define the intertwining operator $A(s, \lambda)$ on  $\mathfrak{H}_\lambda$ by
 \beq\label{Intw}
 A(s, \lambda)F(x)=\int_{V_s}F(v_ssx)\,dv_s.
 \eeq
 It is easy to see formally that $A(s, \lambda)f\in \mathfrak{H}_{s^{-1}\lambda}$ and  for any $g\in G$
 \beq
  A(s, \lambda)\pi_\lambda(g)=\pi_{s^{-1}\lambda}(g)A(s, \lambda).
 \eeq
 In particular, $ A(s, \lambda)\varphi _\lambda=c_s(\lambda)\varphi _{s^{-1}\lambda}$, where
 \beq\label{c-s}
 c_s(\lambda)=\int_{V_s}\phi _\lambda(v_s)\,dv_s=\int_{V_s}e^{\<i\lambda+\rho, H(v_s)\>}\,dv_s.
 \eeq
 Elementary arguments together with the well-known convexity properties of the Cartan component $H(x)$ in the Iwasawa decomposition assure the convergence of all integrals in the tubular domain $\Im \lambda \in C_+$ in $\mathfrak{a}^*_\mathbb{C}$ (where $C_+\subset \mathfrak{a}$ is the positive Weyl chamber).
 The analytic continuation of the integrals \eqref{Intw}, \eqref{c-s} (as   meromorphic functions of $ \lambda$) will be discussed below. When $V_s=V$ (this happens when $s=w_0$ is the longest Weyl group element which maps $\Delta_+$ onto $ \Delta_-,$), $c_s(\lambda)=c(\lambda)$ is the famous \emph{Harish-Chandra function.}

 Let us normalize the intertwining    operators by setting
 \beq\label{Intw-norm}
  B(s, \lambda)f(x)=\frac{1}{c_s(\lambda)}\int_{V_s}f(v_ssx)\,dv_s.
 \eeq
 Thus we have
 \beq\label{sph-s}
 B(s^{-1}, \lambda)\varphi_\lambda=\varphi _{s\lambda}.
 \eeq
  The uniqueness of the  Whittaker vector implies  also that
 \beq\label{wh-s}
 B(s^{-1}, \lambda)w_{\lambda,f}=M(s, \lambda,f)w_{s\lambda,f},
 \eeq
 where the coefficient $M(s, \lambda,\psi)=\overline{M(s, \overline{\lambda},\overline{f})}$ is some meromorphic function. Formulae \eqref{sph-s},  \eqref{wh-s} immediately imply functional equations for spherical and Whittaker functions:
 \beq\label{func_eq}
 \begin{aligned}
 \phi_\lambda(x)&=\phi_{s\lambda}(x), \\
  W_{s\lambda }(x)&=M(s, \lambda,f)^{-1}W_{\lambda   }(x) \q \text{for all}\; s\in W.
 \end {aligned}
 \eeq

 The Harish-Chandra function $c(\lambda)$ plays the central r\^ole in the spectral theory of spherical   functions. In his first paper \cite{HC} Harish-Chandra did not manage to prove the full range of properties of $c(\lambda)$ needed to establish the Plancherel theorem. It came as a total surprise a few years later that the integrals \eqref{c-s} may be calculated explicitly\footnote{One notable exception are rank one groups which have been explored already in the first paper of Harish-Chandra;  here $c(\lambda)$ may be computed in a  routine way using the spectral theory of Sturm--Liouville operators; the \emph{c}-function in this case coincides with the Jost function for the radial  part of the Laplace--Beltrami operator.}. The famous Gindikin--Karpelevich formula, based on an ingenious change of variables in the unipotent group $V$, gives for $c_s(\lambda)$ a simple expression in terms of the product of Euler's Gamma-functions:
 \beq
 \begin{aligned}
 c_s(\lambda)&=\prod_{\alpha\in \Delta(s)}c_\alpha(\lambda), & c_\alpha(\lambda)&=\frac{\Gamma(\lambda_\alpha)\Gamma(^1/_2)}{\Gamma(\lambda_\alpha\,+\, ^1/_2)}, & \lambda_\alpha&=\frac{\<\lambda, \alpha\>}{\<\alpha, \alpha\>}.
\end {aligned}
 \eeq
 In particular,
 \beq
 c(\lambda)=\prod_{\alpha\in \Delta_+}c_\alpha(\lambda).
 \eeq
These formulae provide also an explicit analytic continuation of the Harish-Chandra function (initially defined by an integral which converges only for $\Im\lambda\in C_+$).  The same change of variables allows to establish a factorization formula and an   analytic continuation for the intertwining    operators reducing them to the rank 1 case. Consequently, the functions $M(s, \lambda,\psi)$ which enter the functional equation for the Whittaker functions are also products of functions on one variable. Explicitly we have
 \begin{teo}
 (i) For any $s_1, s_2\in W$ we have
 \beq
 M(s_1 s_2, \lambda,\psi)= M(s_2, \lambda,\psi)M(s_1, s_2\lambda,\psi).
  \eeq
 (ii) For $\alpha\in P$ let $s_\alpha\in W$ be the corresponding elementary reflection. Then
\begin{gather}
 M(s_\alpha, \lambda,f)=e_\alpha(\lambda)e_\alpha(-\lambda)^{-1}\left[\frac{|f (e_\alpha)|}{2\sqrt{2\<\alpha, \alpha\>}}\right]^{2\alpha(\lambda)}, \;\text{where}\\
 e_\alpha(\lambda)=2^{1-\lambda_\alpha}\sqrt{\pi}\,\Gamma({\lambda_\alpha}\,+\,^1/_2).
 \end{gather}
 \end{teo}
 Both the \emph{c}-function and the \emph{M}-function appear in the formula for the scattering matrix for the Toda lattice. The phase factor which compares the behavior of the incoming and the outcoming waves inside the positive Weyl chamber is given by
 \beq
S_{w_0}(\lambda)=\frac{c(\lambda)}{c(w_0\lambda)}M(w_0, \lambda),
 \eeq
 where $w_0$ is the longest element in the Weyl group.  For spherical functions this phase factor is given just by the ratio of \emph{c}-functions,
  \beq
S^0_{w_0}(\lambda)=\frac{c(\lambda)}{c(w_0\lambda)}.
 \eeq

 \begin{remark}The QISM technique which we discuss is the next section naturally leads to a different normalization of Whittaker functions which satisfy a much simpler functional equation
 \beq\label{simpl}
 \psi_{s\lambda }(x)= \psi_{\lambda }(x) \q \text{for all}\; s\in W.
  \eeq
  Since the spectrum multiplicity for the Toda lattice is 1, functions $W_{\lambda   }$ and $\psi_{\lambda }$ should be proportional.
  Remarkably, the symmetry with respect to the Weyl group can be achieved by multiplication by an explicit meromorphic function \cite{Leb1}. Namely, the following simple result holds true:
 \end{remark}
 \begin{pro}
 Let $b(\lambda)$ be the denominator of the Harish-Chandra function,
 $$\label{denom}
 b(\lambda)=\prod_{\alpha\in\Delta_+}\Gamma\left(\frac{\<\alpha, \lambda\>}{i}+\frac{1}{2}\right).
 $$
 Then the modified Whittaker function $\psi_\lambda(x)= b(\lambda)W_{\lambda   }(x)$ satisfies the functional equation \eqref{simpl}.
\end{pro}

 The Gindikin--Karpelevich formula created a sort of new paradigm in Harmonic analysis on semisimple Lie groups. It became a common wisdom that the asymptotical behavior of matrix coefficients of irreducible representations, such as spherical or Whittaker functions, is described  by simple factorized formulae. Similar results hold also for Lie groups over local fields; their extension to Lie groups over the rings of ad\`eles leads to a highly romantic Langlands program (or at least to some part of it, as exposed in the famous book of Gelfand and Pyatetsky-Shapiro \cite{GP}); in this case \emph{c}-functions and \emph{M}-functions are replaced by  zeta or L-functions given by appropriate Euler products and satisfying  the same type of functional equations. On the other hand, the eigenfunctions of the Laplace operators themselves are given by much less manageable integral expressions or infinite sums involving complicated functions of several variables. The rather radical change of this paradigm brought about by the Quantum Inverse Scattering Method consists in the following main points:
\begin{enumerate}
  \item[(i)] There exists a representation (in the sense of Dirac) in which the algebras of quantum integrals of motion for the Toda lattice, or of the Sutherland system are freely generated by first order difference operators.
  \item[(ii)] In this representation the common eigenvectors of these algebras are decomposable, i.e., they are products of functions of one variable.
  \item[(iii)] The transformation operators which relate this representation to the standard coordinate representation may be constructed recursively and amount to a sequence of ordinary Fourier transforms.
\end{enumerate}
Thus in this new representation not only the asymptotics of the Whittaker functions or of the spherical functions are given by product formulae of the Gindikin--Karpelevich type: the Whittaker functions and the spherical functions themselves are given by simple product formulae involving products of Gamma functions. The Gamma functions now appear, for the first time in the context of representation theory, directly as eigenfunctions of difference operators.

As already noted, this new approach to the  Whittaker functions and the quantum Toda lattice  is  largely due to the ITEP group (Lebedev, Kharchev, Gerasimov, Oblezin). Their work was preceded by quite a few very important contributions. An early attempt to understand the quantum Toda lattice is due to  Gutzwiller \cite{Gutzwiller}. A few years  later Sklyanin reached a very profound understanding of   Gutzwiller's approach basing on his seminal technique of quantum separation of variables \cite{Skl1}, \cite{Skl2}. An important contribution is also due to Pasquier and Gaudin \cite{Pasq}. The real challenge which has been addressed in all these papers is the spectral problem for the periodic quantum  Toda lattice. The separation of variables applies in this case as well, but instead of simple difference equations of order one we now get much more complicated equations of order two. Application of the quantum separation of variables to the representation theory, which is technically more elementary,  appears as a rather unexpected but highly interesting byproduct  of this method.

The recursive construction of transformation operators referred to in (iii) is based on the following key idea, again in the spirit of Dirac: the distribution kernel of a unitary transformation operator we are seeking may be regarded as a complete system of eigenfunctions for some auxiliary Hamiltonian. It turns out that this auxiliary Hamiltonian is again a Toda Hamiltonian with one potential term corresponding to the last simple root removed.  Thus the eigenfunctions of this Hamiltonian are just the Toda eigenfunctions in rank $n-1$ multiplied by a free exponential associated with the removed root. Iterating this procedure we can express  the Toda eigenfunctions via a sequence of ordinary Fourier transforms of an aggregate composed of Gamma functions. (Importantly, the number of integrations exceeds the number of independent variables!) The key technical point which sets this machinery into motion is the explicit formula giving the action of the rank $n$ Toda Hamiltonian on the eigenfunctions in rank $n-1$. It is this point that represented the key initial contribution of QISM to the problem in question.

 In order to construct this  representation  we shall give an alternative description of the quantum Toda model based on QISM. The key result described in section \ref{QISM} is the explicit expression for the wave functions of the open Toda lattice in the separation representation (see formulae \eqref{wf_red}, \eqref{wf} below)  together with  the   formulae \eqref{MB}, \eqref{MB-long} which give the wave functions of the Toda lattice in the standard coordinate representation.   In    section \ref{Ug} we shall describe the associated representations of the universal enveloping algebra.

 \section{Quantum Toda Lattice:  the point of view of QISM}\label{QISM}

The application of the Quantum Inverse Scattering Method to the Toda lattice is based on the use of a $2\times2$ Lax pair which goes back to Flashka and Manakov \cite{Fl}, \cite{Man}\footnote{The use of this Lax pair imposes a considerable restriction of generality: we are now dealing with the series $A$ root system and hence with the case of $SL(N)$ or $GL(N)$.  The generalization to other classical series is also possible, although it requires some extra work; working in full generality, i.\,e., in a classification independent way for arbitrary semisimple Lie groups requires probably some new ideas. }. Its quantum version was discovered in the early days of QISM; the main problem discussed at that time was the spectrum of \emph{periodic} Toda lattice with the Hamiltonian
\beq\label{T_perio}
H=-\Delta+\sum_{k=1}^{N-1}e^{x_{k+1}-x_{k}} + e^{x_{1}-x_{N}}, \q x\in \mathbb{R}^N, \q \sum_k x_k=0.
\eeq
which differs from \eqref{T_open} by the extra term in the potential. Since the potential energy in this case is growing in all directions, the spectrum of the {periodic} Toda lattice is discrete. Its explicit calculation proved to be much more complicated than in other models, such as the Heisenberg magnet, which were resolved at that time using the algebraic Bethe ansatz technique. From the point of view of Representation theory the periodic Toda model differs rather dramatically from the open one: in this case the underlying Lie group is infinite dimensional (it is the loop group of $SL(N)$). Technically, this means that the underlying Lax pair contains a \emph{spectral parameter}. Spectral parameter plays of course the key role in the QISM approach as well. The algebraic treatment of the Toda model in QISM formalism starts with the rather \emph{ad hoc}
 choice of the Lax matrix. Let $q_m, p_m,   \, m=1, \dots  n, $ be canonical position and momentum operators  with commutation relations\footnote{In various papers on QISM (cf. \cite{Skl1} and \cite{Leb1}) ) the Planck constant is usually retained as an extra scaling parameter. To simplify the formulae and also to facilitate the comparison with representation theory throughout this paper we assume that $\hbar=1$, both in the Toda Hamiltonian and in the quantum R-matrix.}
 $$
 [p_j, q_k]=-i\delta_{jk}, \; i=\sqrt{-1}.
 $$
 We introduce a family of $2\times2$-matrices with operator entries depending on a complex parameter $u$,
 \beq
 L_m(u)=\begin{pmatrix}
 u-p_m&-e^{q_m}\\
 e^{-q_m}&0
 \end{pmatrix}.
 \eeq
 We regard $ L_m(u)$ as an element of $\mathrm{Mat}(2, \mathbb{C})\otimes \mathrm{End} (\mathfrak{H})$, where $\mathfrak{H}$ is the representation space for the canonical operators $q_m, p_m$. The space $\mathfrak{H}$ is commonly referred to as the quantum space and the space $\mathbb{C}^2$ as the auxiliary space. The success of the QISM  machinery relies on a clever way to encode the commutation relations in the quantum space by an appropriate linear transformation in the  auxiliary space. Let $R(u)$ be a linear operator acting in $\mathbb{C}^2\otimes\mathbb{C}^2$,
 \beq
 R(u)=uI-iP,
 \eeq
where $I$ is the identity operator and $P\in \mathrm{End}(\mathbb{C}^2\otimes\mathbb{C}^2)$ is the permutation operator, $P(v\otimes w)=w\otimes v$.
\begin{pro}
(i) We have
\beq\label{CCR}
\left[p_m, e^{\pm q_k}\right]=\mp ie^{\pm q_k}\delta_{mk}, \q i=\sqrt{-1}.
\eeq
(ii) Put $L^1_m(u)=L_m(u)\otimes I, \; L^2_m(u)=I\otimes L_m(u)$; we regard  $L^1_m(u),  L^2_m(u)$ as operators in $\mathrm{End}(\mathbb{C}^2\otimes\mathbb{C}^2)\otimes\mathrm{End} (\mathfrak{H})$. The commutation relations \eqref{CCR} are equivalent to the operator identities
\beq\label{RLL}
\begin{aligned}
L^2_m(v)L^1_m(u)&=R(u-v)L^1_m(u)L^2_m(v)R(u-v)^{-1},  \\
 L^1_m(u)L^2_k(v)&=L^2_k(v)L^1_m(u) \; \text{for} \; m\neq k.
 \end{aligned}
\eeq
(iii) Let $T_N(u)=L_N(u)L_{N-1}(u)\dots L_1(u)$ be the monodromy matrix associated with `local'  Lax matrices $L_m(u)$; we put   as usual $T^1(u)=T(u)\otimes I$, $T^2(u)=I\otimes T(u)$. The monodromy matrices satisfy
\beq\label{Hopf}
\begin{aligned}
T_N^2(v)T_N^1(u)&=R(u-v)T_N^1 (u)T_N^2 (v)R(u-v)^{-1},
 \end{aligned}
\eeq
(iv) We write
\beq
\begin{aligned}
T_N(u)&=\begin{pmatrix}
A_N(u)&B_N(u)\\
C_N(u)&D_N(u)
\end{pmatrix}, & t(u)&= \mathrm{tr }\,T_N(u)=A_N(u)+D_N(u).
 \end{aligned}
\eeq
Commutation relations \eqref{Hopf} imply that
\beq
[t(u), t(v)]=0.
\eeq
\end{pro}
 The commutative family of operators $t(u)$ is clearly a polynomial in $u$,
 $$
 t(u)=u^N+t_1u^{N-1}+t_2u^{N-2}+\dots,
 $$
 where in particular
\beq
\begin{aligned}
 t_1&=-P_N=-p_1-\dots+p_N, \\ t_2&=\frac12 P_N^2-\frac12(p^2_1+\dots+p^2_N)-e^{x_2-x_1}-\dots e^{x_N-x_{N-1}}-e^{x_1-x_N}.
 \end{aligned}
\eeq
Thus the coefficient $t_1$ is (up to sign) the full momentum and $t_2$ gives essentially the Hamiltonian of the periodic Toda lattice of length $N$.

For applications to open Toda lattice (and hence eventually to Representation  theory) we need other coefficients of the monodromy matrix. Commutation relations \eqref{Hopf} imply that
\beq
\begin{aligned}\label{CommAC}
 \left[B_N(u),B_N(v)\right]&=0, \;  \left[C_N(u),C_N(v)\right] =0,\\
(u-v+i)A_N(v)C_N(u)& =(u-v)C_N(u)A_N(v)+i A_N(u)C_N(v).
\end{aligned}
\eeq
It is easy to see that the coefficient $A_N(u)$ is the generating function for the Hamiltonian of the open Toda lattice of length $N$ and of its integrals of motion.  We can write
\beq\label{AD}
\begin{aligned}
 A_{N}(u)&=u^N+\sum_{m=1}^{N}X_m u^{N-m}, &
D_{N}(u)&=\sum_{m=2}^{N}Y_m u^{N-m},
\end{aligned}
\eeq
where $X_m, Y_m$ are two sets of mutually commuting operator coefficients. The coefficients $X_1=-P$ and $X_2=\frac12P^2-H$ give  the total momentum and the Hamiltonian for the open Toda lattice; the coefficients $Y_m$ add  up to $X_m$ to yield the  quantum integrals of motion  for the periodic Toda.
Moreover,
\beq\label{AC}
\begin{aligned}
A_N(u)&=(u-p_N)A_{N-1}(u)+e^{-x_N}C_{N-1}(u),\\
C_N(u)&=e^{-q_N}A_{N-1}(u)
\end{aligned}
\eeq
Note that $A_{N-1}(u)$, which is the generating function for the Hamiltonians of the open Toda lattice of length $N-1$, commutes with $q_N, p_N$. Clearly, $C_N(u)$ is a polynomial in $u$ and hence it may written in factorized form,
\beq\label{roots}
C_N(u)=e^{-q_N}\prod_{m=1}^{N-1}\left(u-\hat{\lambda}_m\right),
\eeq
where $\hat{\lambda}_m$, $m=1, \dots, N-1,$ are mutually commuting self-adjoint operators (``the operator-valued roots of $C_N(u)$'').
\begin{definition}\label{Sep}Operators $\hat{\lambda}_m$ are called\emph{ quantum separated variables}; spectral representation for these operators is called \emph{separation representation.}
\end{definition}
In order to justify this definition we need to introduce some operator calculus.
Following \cite{Skl1}, we define operators $\Lambda_j^\pm$, $j=1, \dots, N-1$ be setting
\beq\label{lambda}
\begin{aligned}
\Lambda_j^-&=\left.A_{N}(u)\right|_{u=\hat{\lambda}_j}=\hat{\lambda}_j^N+\sum_{m=1}^{N-2}X_m \hat{\lambda}_j^{N-2-m},\\
\Lambda_j^+&=\left.D_{N}(u)\right|_{u=\hat{\lambda}_j}=\sum_{m=2}^{N}Y_m \hat{\lambda}_j^{N-m},
\end{aligned}
\eeq
where $X_m, Y_m$ are  defined in \eqref{AD}.  (Note that  $X_m, Y_m$ do not necessarily commute with $\hat{\lambda}_j$; the definition \eqref{lambda} fixes the `left ordering' of the  coefficients.) The coefficients of the monodromy matrix may be restored  from $\hat{\lambda}_j, \Lambda_j^\pm$ by means of the Lagrange interpolation formulae,
\begin{equation}\label{Lagr}
\begin{aligned}
A_N(u)&=\left(u-P+\hat{\lambda}_1+\dots \hat{\lambda}_{N-1}\right)\prod_{j=1}^{N-1}(u-\hat{\lambda}_j)+
\sum_{j=1}^{N-1}\left(\prod_{\begin{smallmatrix}
k=1\\
k\neq j
\end{smallmatrix}}^{N-1}\frac{u-\hat{\lambda}_k}{\hat{\lambda}_j-\hat{\lambda}_k}\Lambda_j^-\right),\\
D_N(u)&=\sum_{j=1}^{N-1}\left(\prod_{\begin{smallmatrix}
k=1\\
k\neq j
\end{smallmatrix}}^{N-1}\frac{u-\hat{\lambda}_k}{\hat{\lambda}_j-\hat{\lambda}_k}\Lambda_j^+\right).
\end{aligned}
\end{equation}
The basic commutation relations \eqref{Hopf}, \eqref{CommAC} imply that
\beq\label{comm}
\begin{aligned}
\left[\hat{\lambda}_j,\, \hat{\lambda}_k\right]&=\left[\Lambda_j^\pm, \,\Lambda_k^\pm\right]=0,\\
\Lambda_j^\pm  \hat{\lambda}_k&=\left(\hat{\lambda}_k\pm i\delta_{jk}\right)\Lambda_j^\pm, \q i=\sqrt{-1}.
\end{aligned}
\eeq
These key commutation relations mean  that in the spectral representation for $\hat{\lambda}_k$ the operators $\Lambda_j^\pm$ act as translation operators affecting only one variable, which justifies definition \ref{Sep}. Explicitly, this spectral representation is given by
\beq
\begin{aligned}\label{spec}
\hat{\lambda}_j f (p, \lambda_1, \lambda_2, \dots, \lambda_{N-1})&=\lambda_j\cdot f(k, \lambda_1, \lambda_2, \dots, \lambda_{N-1}),\\
\Lambda_j^\pm f(p, \lambda_1, \lambda_2, \dots, \lambda_{N-1})&=i^{\pm N}f(k, \lambda_1, \lambda_2, \dots, \lambda_j\pm i, \dots, \lambda_{N-1});
\end{aligned}
\eeq
where we have added one more variable $p$, the eigenvalue of the total momentum $P$. The phase factor $i^{\pm N}$ is introduced for convenience.  Together with the conjugate variable $x_N$ the momentum  $P$ completes the set of independent variables to yield a representation of the  rank $N$ Heisenberg algebra. The action of $P, x_N$ in the spectral representation is given by
\beq
\begin{aligned}
P f (p, \lambda_1, \lambda_2, \dots, \lambda_{N-1})&=pf (p, \lambda_1, \lambda_2, \dots, \lambda_{N-1}),\\
e^{-x_N} f(p, \lambda_1, \lambda_2, \dots, \lambda_{N-1})&=f (p-i, \lambda_1, \lambda_2, \dots, \lambda_{N-1}).
\end{aligned}
\eeq
\begin{remark} The roots of the polynomial equation \eqref{roots} are defined only up to a permutation; it is thus natural to assume that functions $f$ in the spectral representation \eqref{spec} are symmetric with respect to the variables $\lambda_j$. With this symmetry condition imposed the Toda wave functions will automatically satisfy the simple functional equation \eqref{simpl} (cf. remark \ref{func_eq} above).
\end{remark}

The use of separation representation provides a rather novel point of view on the computation of the Plancherel measure. The standard way to compute the Plancherel measure is based on the study of the coordinate asymptotics of the wave packets \eqref{wpack}. In this way one gets the famous Harish-Chandra formula for the Plancherel measure which is valid both for spherical functions and   for the Whittaker functions (normalized as in \eqref{whit}):
\beq
\mu_0(\lambda)\,d\lambda=\frac{1}{|c(\lambda)|^2}\,d\lambda.
\eeq
Factorization properties of \emph{c}-function as well as the difference equation for Gamma functions do not play any r\^ole in the  calculation of this measure. As noted by Sklyanin \cite{Skl1}, the use of the separation representation  allows to compute the Plancherel measure directly using the conjugation properties of the operators $\hat{\lambda}_j, \Lambda_j^\pm$. Let us denote by $T_N^*(u)$ the monodromy matrix obtained  from $T_N(u)$ by hermitian conjugation in the quantum space (which acts on matrix coefficients of $T_N(u)$ only). We have
\beq\label{inv}
\begin{aligned}
T_N^*(u)&=T_N(\bar{u}),
\end{aligned}
\eeq
 which immediately yields
 \beq\label{conj}
\begin{aligned}
\hat{\lambda}_k^*&=\hat{\lambda}_k,\\
\left(\Lambda_j^\pm\right)^*&=\prod_{\begin{smallmatrix}
k=1\\
k\neq j
\end{smallmatrix}}^{N-1}\frac{\hat{\lambda}_j-\hat{\lambda}_k\pm i}{\hat{\lambda}_j-\hat{\lambda}_k}\Lambda_j^\pm.
\end{aligned}
\eeq
In order to reconcile these conjugation properties with the spectral representation \eqref{spec} we assume that the representation space is equipped with the inner product of the form
\begin{equation*}
\begin{aligned}
\<f,g\>&=\int_{\mathbb{R}^N} f (p, \lambda_1,   \dots, \lambda_{N-1})  \overline{g(p, \lambda_1,  \dots, \lambda_{N-1})}\mu(\lambda_1,   \dots, \lambda_{N-1})\,dp  d\lambda,
\end{aligned}
\end{equation*}
where $\mu(\lambda_1,   \dots, \lambda_{N-1})$ is some positive symmetric function. Now a comparison with \eqref{conj} yields  a system of difference equations for $\mu$:
\beq
\begin{aligned}
\mu(\lambda_1,   \dots, \lambda_j+i, \dots \lambda_{N-1})&=\prod_{\begin{smallmatrix}
k=1\\
k\neq j
\end{smallmatrix}}^{N-1}\frac{{\lambda}_j-{\lambda}_k- i}{{\lambda}_j-{\lambda}_k}, \q j=1, \dots, N-1,
\end{aligned}
\eeq
which has a solution
\beq\label{nom}
\begin{aligned}
\mu(\lambda_1,   \dots,   \lambda_{N-1})&=\prod_{\begin{smallmatrix}
k=1\\
j<k
\end{smallmatrix}}^{N-1}\frac{1}{\left|\Gamma\left(\frac{{\lambda}_j-{\lambda}_k}{i}\right)\right|^2}.
\end{aligned}
\eeq
The denominator in \eqref{nom} is precisely the nominator of the Harish-Chandra function; formula \eqref{nom} has to be compared with the factor in \eqref{denom} which changes the normalization of the Whittaker function in order to make it symmetric. Thus the symmetry with respect to involution \eqref{inv} determines the Plancherel measure completely (up to a numerical factor).

Interpolation formulae \eqref{Lagr} imply that in the separation representation the simultaneous diagonalization of the quantum integrals of motion  for the open Toda lattice  is reduced to a system   of difference equations. We can parameterize the common eigenfunctions $\psi_\alpha$,  $\alpha=(\alpha_1, \dots, \alpha_N)$ of these integrals by assuming that the eigenvalues of the operators $X_m$ introduced in  \eqref{AD} are elementary symmetric polynomials,
\beq
X_m\psi_\alpha=(-1)^m\sigma_m(\alpha_1, \dots, \alpha_N)\psi_\alpha
\eeq
This parametrization   clearly agrees with the Harish-Chandra homomorphism which we used to describe quantum integrals of motion in section~\ref{rep}. Equivalently this means that
\beq
A_N(u)\psi_\alpha=\prod_{j=1}^N(u-\alpha_j)\cdot \psi_\alpha.
\eeq
According to \eqref{Lagr}, the commutative pencil $A(u)$ is completely determined by its values $A(\lambda_j)=\Lambda_j^-$, $ j=1,\dots, N-1$ together with the operator of the total momentum. Since in the separation representation the total momentum acts as multiplication operator,   the common eigenfunctions of $A(u)$ have the form
\beq\label{wf}
\psi_\alpha(P, \lambda_1, \dots, \lambda_{N-1})=\delta(P-\alpha_1)\phi_\alpha( \lambda_1, \dots, \lambda_{N-1}),
\eeq
where $\phi_\alpha$ satisfies the system of difference equations
\beq\label{dif}
\begin{aligned}
\Lambda_j^-\phi_\alpha( \lambda_1, \dots, \lambda_{N-1})&=\prod_{k=1}^N(\lambda_j-\alpha_k)\phi_\alpha( \lambda_1, \dots, \lambda_{N-1}), \; j=1, \dots, N-1.
\end{aligned}
\eeq
Its obvious solution is
\beq\label{wf_red}
\begin{aligned}
\phi_\alpha( \lambda_1, \dots, \lambda_{N-1})&=\prod_{j=1}^{N-1}\prod_{k=1}^N \Gamma\left(\frac{\lambda_j-\alpha_k}{i}\right).
\end{aligned}
\eeq
Formulae \eqref{wf_red}, \eqref{wf} give an explicit expression for the wave functions of the open Toda lattice in the separation representation.
\begin{remark} Quantum separation of variables works also for the periodic Toda lattice; in fact, it is in this context that it has been originally conceived \cite{Skl1}. The generating function for the integrals of motion of the periodic Toda lattice is $t_N(u)=A_N(u)+D_N(u)$. Again, the total momentum acts trivially and may be excluded. Parameterizing   the common eigenfunctions of quantum integrals of motion (in the center-of-mass frame) by the roots of  $t_N(u)$ (now denoted by $E=(E_1,\dots, E_N)$,  to distinguish them from the former case)  we get for the eigenfunctions of the periodic Toda the following system of compatible difference equations:
\beq\label{periodic}
\begin{aligned}
\left(\Lambda_j^+ +\Lambda_j^-\right)\Phi_E(\lambda_1, \dots, \lambda_{N-1})&=\prod_{k=1}^N(\lambda_j-E_k)\Phi_E( \lambda_1, \dots, \lambda_{N-1}), \\ j&=1, \dots, N-1.
\end{aligned}
\eeq
In contrast with \eqref{dif} these equations are of order 2, which makes their solution much less elementary\footnote{Important results in connection with this difference system may be found in \cite{Pasq}, \cite{Leb2}}. Recall that in the classical case the Lax pair for the periodic Toda is provided by  infinite-dimensional loop algebras and  its solution  requires the machinery of Algebraic Geometry; the case of the open Toda lattice corresponds to the degeneration   of the spectral curve associated with the Lax pair into a singular rational curve; in this degenerate case its solution is reduced to matrix factorization. In the quantum case a similar degeneration consists in the replacement of the second order system \eqref{periodic} by a much simpler first order system  \eqref{dif}. Very remarkably, this simple  system yields plenty of information and contains much of the Harish-Chandra theory.
\end{remark}
In order to complete the solution  of the open Toda lattice we must describe the transformation operator  which relates the separation representation to the standard coordinate representation. Its description is already implicit in formula \eqref{AC} which states that the diagonalization of the operators $C_N(u)$ is equivalent to the solution of the rank $N-1$ open Toda lattice (up the introduction of an extra coordinate variable).
 \begin{pro}
 Let $\psi_\lambda(x_1,\dots,x_{N-1})$, $\lambda=(\lambda_1, \dots, \lambda_{N-1})$, be the common eigenfunction of the operator pencil $A_{N-1}(u)$. Put
 \beq\label{Psi}
 \Psi_{\lambda, P}(x_1,\dots,x_{N})=\exp\left( iP-i\sum_{m=1}^{N-1}\lambda_m\right)x_N \cdot \psi_\lambda(x_1,\dots,x_{N-1}).
 \eeq
 We have
 \beq
 C_N(u)\Psi_{\lambda, P}=-e^{x_N}\prod_{j=1}^{N-1}(u-\lambda_j)\Psi_{\lambda,P}.
 \eeq
\end{pro}
The phase factor in \eqref{Psi} is chosen in such a way that $\Psi_{\lambda, P}$ is an eigenfunction of the operator of total momentum with eigenvalue $P$.

Speaking more accurately, the wave functions which correspond to the operators with  continuous spectrum (of course, such wave functions do not lie in the Hilbert space) should be interpreted as distribution kernels of transformation operators which rely different sets of quantum observables. Following Dirac, we shall label different ``coordinate systems''  in the space of quantum states (alias, quantum representations) by indicating the complete sets of observables which are diagonal, i.e., act by multiplication. We are dealing with the coordinate representation in which the diagonal operators are the standard coordinates, the separation representation in which the diagonal operators are the operator roots $\lambda=(\lambda_1, \dots, \lambda_{N-1})$ of $C_N(u)$ of the total momentum, and with the spectral representation, in which the diagonal operators are the   operator roots  $\alpha=(\alpha_1, \dots, \alpha_{N})$ of $A_N(u)$. In Dirac notation the wave function which diagonalizes the open Toda Hamiltonians may be written simply as $\<x|\lambda\>$
$$
\psi_\alpha(x_1, \dots, x_N)=\<x|\alpha\>.
$$
In a similar way, the notation for $\Psi_{\lambda, \lambda_N}(x_1,\dots,x_{N})$ may be condensed to
$$
 \Psi_{\lambda, \lambda_N}(x_1,\dots,x_{N})=\<x|\lambda\>
$$
with $\lambda=(\lambda_1, \dots, \lambda_{N})$, as in \eqref{Psi}. Finally, the wave function \eqref{wf} in separation representation is condensed to
$$
\psi_\alpha(P, \lambda_1, \dots, \lambda_{N-1})=\<\lambda,P|\alpha\>.
$$
This interpretation immediately yields  the following integral formula for the wave function of the open Toda lattice in the coordinate representation:
\begin{multline}\label{MB}
\psi_\alpha(x_1, \dots, x_N):=\<x|\alpha\>=\int \<x|\lambda\>\<\lambda|\alpha\>\, \mu(\lambda)\,d\lambda\,dp=\\
\int e^{\left( iP-i\sum_{m=1}^{N-1}\lambda_m\right)x_N}\prod_{j=1}^{n-1}\prod_{k=1}^N \Gamma\left(\frac{\lambda_j-\alpha_k}{i}\right) \cdot \psi_\lambda(x_1,\dots,x_{N-1})\mu(\lambda)\,d\lambda\,dP,
\end{multline}
where the measure density $\mu(\lambda)$ is given by \eqref{nom}.

Formula \eqref{MB} immediately leads to a simple recursion which ends up with an explicit formula for $\psi_\alpha(x_1, \dots, x_N)$ in terms of multiple Fourier transforms of an agregate composed of Gamma functions. This is the famous Mellin--Barnes type formula\footnote{The classical Mellin--Barnes   formula gives an integral representation for the Whittaker function of one variable which corresponds to the case of $\mathfrak{sl}(2)$. } first obtained by Kharchev and Lebedev in \cite{Leb1}. In order to write down this formula  explicitly let us rearrange the integration variables into a triangular array
$ \lambda_n= (\lambda_{n1}, \lambda_{n2}, \dots, \lambda_{nn})$,    $n=1, \dots N$; the upper raw $(\lambda_{N1}, \lambda_{N2}, \dots, \lambda_{NN})$ is identified with our spectral variables $\alpha=(\alpha_1, \dots, \alpha_{N})$. With this notation we get
\begin{multline}\label{MB-long}
\psi_{\lambda_N}(x_1, \dots, x_N)=\\
\int\prod\limits_{n=1}^{N-1}\frac{\prod\limits_{k=1}^n\prod\limits_{m=1}^{n+1}\Gamma\left(\frac{\lambda_{nk}-\lambda_{n+1, m}}{i}\right)}{\prod_{\begin{smallmatrix}
s\neq p,   \\s,p\leq n
    \end{smallmatrix}}\Gamma\left(\frac{\lambda_{ns}-\lambda_{n, p}}{i}\right)}\exp\left\{i\sum_{n=1}^N\sum_{k=1}^{n-1}(\lambda_{nk}-\lambda_{n-1,k})\right\}\prod\limits_{\begin{smallmatrix}
    n=1\\j\leq n
 \end{smallmatrix}   }^{N-1}d\lambda_{nj}.
\end{multline}
Formula \eqref{MB-long} represents the main result of the QISM technique as applied to the open Toda lattice. Integration in \eqref{MB-long} is performed over an appropriate contour which avoids singularities of Gamma functions; the explicit character of this formula makes the convergence control quite simple. The study of  the asymptotic behaviour  of Whittaker functions becomes elementary too and is reduced to the Cauchy residue theorem.

 \section{Whittaker vectors  and spherical vectors in Gelfand--Zetlin representation}\label{Ug}

Quantum separation of variables described in the previous section has got a profound and unexpected connection with a classical construction in Representation theory which goes back to an idea of Gelfand and Zetlin \cite{GZ}.  Consider a chain of natural embeddings of matrix algebras
$$\mathfrak{gl}(N)\supset \mathfrak{gl}(N-1) \supset\dots \supset \mathfrak{gl}(1),$$
 where an $(n-1)\times (n-1)$-matrix is identified with the upper left corner of an $n\times n$-matrix. The centers $\mathcal{Z}_n\subset U(\mathfrak{gl}(n))$, $n=N, N-1, \dots 1$, of the embedded universal enveloping algebras generate a maximal abelian subalgebra   $\mathcal{A}_{GZ}\subset U(\mathfrak{gl}(N))$. For each $n$ the center $\mathcal{Z}_n\subset U(\mathfrak{gl}(n))$ is isomorphic to the algebra of symmetric polynomials in $n$ variables; thus it is natural to parameterize the spectrum of  $\mathcal{A}_{GZ}$ by a triangular array of variables $ \lambda_n= (\lambda_{n1}, \lambda_{n2}, \dots, \lambda_{nn})$,    $n=1, \dots N$; on each level the variables  $(\lambda_{n1}, \lambda_{n2}, \dots, \lambda_{nn})$ are defined up to a permutation\footnote{Speaking more accurately, our triangular arrays parameterize the spectrum of a huge algebraic extension of $\mathcal{A}_{GZ}$ with the Galois group $S_1\times S_2\times\dots\times S_N$.}. Arrays of this type have already appeared in formula \eqref{MB-long} which represents the main result of the previous section. We may regard the  variables $ \lambda_n$ as generalized momenta; adding to them the set of conjugate coordinates we get a big Heisenberg algebra\footnote{The elements of the center  $\mathcal{Z}_N\subset U(\mathfrak{gl}(N))$ have no conjugates, so their common spectrum, parameterized by the first line $(\lambda_{N1},\dots, \lambda_{NN})$ of our array,  gives a set of moduli for the entire construction. }. According to an  important observation due to  Gelfand and Kirillov \cite{GK}, the universal enveloping algebra of this Heisenberg algebra is ``almost'' isomorphic to the field of fractions of $U(\mathfrak{gl}(N))$\footnote{This assertion is known as the Gelfand-- Kirillov hypothesis. For real split semisimple Lie algebras it has been proved by  Gelfand and Kirillov in \cite{GK1}.}. Taken seriously, this old algebraic idea leads to a totally new analytic realization of $U(\mathfrak{gl}(N))$ \cite{Leb3}\footnote{On the semiclassical level (with Poisson brackets replacing the commutation relations in $U(\mathfrak{gl}(N))$)  the link between the Gelfand--Zetlin construction and the Gelfand--Kirillov  isomorphism is also discussed in \cite{KW}.}.  Actually, this realization makes use not of the Heisenberg algebra itself which is represented by multiplications and derivations, but rather of the associated Weyl algebra generated by multiplications and translations.

Let $\mathcal{M}$ be the space of   meromorphic functions of  $(\lambda_1, \lambda_{2},  \dots, \lambda_{N-1})$,  $\lambda_{n}\in \mathbb{C}^n$. We define operators $T_{{nj},x}$, $x\in \mathbb{C}$,  acting in the space $\mathcal{M}$ of meromorphic functions of  $(\lambda_1, \lambda_{2},  \dots, \lambda_{N-1})$ by translations,
 $$
 T_{{nj}, x}\Phi(\lambda_1, \lambda_{2},  \dots, \lambda_{N-1})=\Phi(\lambda_1, \dots,    \lambda_{nj}+x, \dots , \lambda_{N-1}).
 $$
 We define a linear representation of $U(\mathfrak{gl}(N))$ in $\mathcal{M}$  which will be called the \emph{Gelfand--Zetlin representation}. The representation operators which correspond to elementary matrices $E_{nm}$ are defined as linear combination with rational coefficients of   translation operators acting in purely imaginary directions.
\begin{theorem}\cite{Leb3} Operators
\beq\label{GZ}
\begin{aligned}
E_{nn}&=\frac{1}{i}\sum_{j=1}^n\lambda_{nj}-\frac{1}{i}\sum_{j=1}^{n-1}\lambda_{n-1,j}, \; n=1,\dots, N-1,\\
E_{n,n+1}&=\frac{1}{i}\sum_{j=1}^n\frac{\prod_{r=1}^{n+1}\left(\lambda_{nj}-\lambda_{n+1,r}- \frac{i}{2}\right)}{\prod_{s\neq j} (\lambda_{nj}-\lambda_{ns})}T_{{nj},-i}, \\
E_{n,n+1}&=\frac{1}{i}\sum_{j=1}^n\frac{\prod_{r=1}^{n-1}\left(\lambda_{nj}-\lambda_{n-1,r}+ \frac{i}{2}\right)}{\prod_{s\neq j} (\lambda_{nj}-\lambda_{ns})}T_{{nj}, i}, \; i=\sqrt{-1},
\end{aligned}
\eeq
define a representation of $U(\mathfrak{gl}(N))$ in $\mathcal{M}$.
\end{theorem}
To prove the theorem one has to check the standard commutation relations between the generators,
\beq\label{rel}
\begin{aligned}
&\left[E_{nn},\, E_{m,m+1}\right]=\left(\delta_{nm}-\delta_{n,m+1}\right)E_{m,m+1},\\
&\left[E_{nn},\, E_{m+1,m}\right]=-\left(\delta_{nm}-\delta_{n,m+1}\right)E_{m+1,m},\\
&\left[E_{n,n+1},\, E_{m+1,m}\right]=\delta_{nm}\left(E_{nn}-E_{n+1,n+1}\right),
\end{aligned}
\eeq
as well as the Serre relations between the generators $E_{n,n+1}$, $E_{m,m+1}$ (which correspond to simple roots of $\mathfrak{gl}(N)$).

As already noted, the first row of the triangular array $$ \lambda_N= (\lambda_{N1}, \lambda_{N 2}, \dots, \lambda_{NN}),$$  which fixes the eigenvalues of the Casimirs in $U(\mathfrak{g}l(N))$, provides a set of moduli for the representation \eqref{GZ}; the same set  of parameters was used in Section \ref{rep} to parameterize the principle series representations.

\begin{remark}
 Representation \eqref{GZ} looks rather remarkable; this is indeed    a totally new type of representations in which elements  of a Lie algebra are acting as functional difference operators rather than as differential operators.  The emergence of difference operators instead of  differential operators is   much more familiar for q-deformed objects such as quantum groups. As a matter of fact, representation \eqref{GZ} is closely related to a representation of an infinite dimensional quantum algebra, the Yangian $Y(\mathfrak{gl}(N))$. We shall return to this question at the end of this section. 
  \end{remark}
Although the   check of commutation relations \eqref{rel} is straightforward, the analytical aspects of this definition (which is motivated by the constructions of QISM) are far from being well-understood.
Moreover, since  the translation operators are acting in  imaginary directions, these difference operators are actually unbounded. The underlying linear space $\mathcal{M}$, which  is chosen \emph{ad hoc}, is tremendously big. However, all these disadvantages and obscure points  are made up by the following wonderful property: Whittaker vectors  and spherical vectors in the Gelfand--Zetlin representation may be chosen in  factorized form, i.e., they are products of meromorphic functions of one variable. More precisely,  the following assertion holds:
\begin{pro} \cite{Leb3} The finite difference equations
\beq\label{Wh_v}
\begin{aligned}
E_{n+1,n}w'_N&=-iw'_N,\\
E_{n,n+1} w_N&=-iw_N, \;  n=1, \dots N-1,
\end{aligned}
\eeq
 admit solutions in  $\mathcal{M}$,
 \beq\label{Wh_sol}
\begin{aligned}
 w'_N&=1,\\
 w_N&=\prod_{n=1}^{N-1} e^{-\pi(n-1)\sum_{j=1}^n \lambda_{nj}}s_n(\lambda_n, \lambda_{n+1}),
\end{aligned}
\eeq
where
\beq
s_n(\lambda_n, \lambda_{n+1})=\prod_{k=1}^n\prod_{m=1}^{n+1}\Gamma\left(\frac{\lambda_{nj}-\lambda_{n+1,m}}{i}+\frac{1}{2}\right).
\eeq
\end{pro}
Clearly, $w_N$, $w'_N$ are Whittaker vectors in the Gelfand--Zetlin representation which correspond to the opposite nilpotent subalgebras.
The proof, which is a bit tricky, uses some identities for symmetric functions (in addition to the standard properties of Gamma).
 \begin{remark}
 Of course, the solutions of difference equations in $\mathcal{M}$ are not unique because of the possibility to multiply a   particular solution by   quasi-constants, which in the present context are  products of periodic functions with purely imaginary periods. It is probably possible to get rid of this freedom by imposing some growth conditions on the solutions, but the question remains unexplored.
 \end{remark}

 The cyclic submodules   $\mathcal{W}, \mathcal{W}'\subset \mathcal{M}$ generated by the Whittaker vectors $w_N$, $w'_N$ are irreducible.
 \begin{pro} (i) The submodules $$\mathcal{W}=U(\mathfrak{gl}(N))w_N, \mathcal{W}'=U(\mathfrak{gl}(N))w'_N$$ are freely generated by the action of the maximal abelian subalgebra $\mathcal{A}_{GZ}\subset U(\mathfrak{gl}(N))$.

 (ii) There exists a natural Hermitian pairing $\mathcal{W}\times \mathcal{W}'\to \mathbb{C}$ which sets the Whittaker modules   $\mathcal{W}, \mathcal{W}$ into duality, so that
 $$
 \<\phi, X\psi\>=- \<X\phi, \psi\> \q \text{for all} \q \phi\in \mathcal{W},\,\psi\in \mathcal{W}', \, X\in \mathfrak{gl}(N)).
 $$
\end{pro}
Explicitly, we have
\beq\label{int}
\<\phi, \psi\>=\int\phi(\lambda) \overline{\psi(\lambda)} \mu(\lambda)\prod_{n=1}^{N-1}\prod_{j=1}^n d\lambda_{nj},
\eeq
where the measure density satisfies the difference equations
\beq
\left(T_{{nj}, i}\mu\right)(\lambda)=\mu(\lambda)\prod_{s\neq j}\frac{\lambda_{nj}-\lambda_{ns}+i}{\lambda_{nj}-\lambda_{ns}}
\eeq
and  is explicitly given by
\beq
\mu(\lambda)=\prod_{n=1}^{N-1}\prod_{s<p}\left(\lambda_{ns}-\lambda_{np}\right)\left(e^{2\pi\lambda_{np}}-e^{2\pi\lambda_{ns}}\right).
\eeq
The integration in \eqref{int} is over the real subspace $\mathbb{R}^{N(N-1)/2}\subset \mathbb{C}^{N(N-1)/2}$; using the standard estimates for the Gamma function one can show that the integral is absolutely convergent  for any $\phi\in \mathcal{A}_{GZ}\cdot w_N$, $\psi\in\in \mathcal{A}_{GZ}\cdot w'_N$.

A similar factorized formula holds for the spherical vectors in  $\mathcal{W}, \mathcal{W}'$.
 \begin{pro} (i) The maximal compact subalgebra $\mathfrak{k}\subset \mathfrak{gl}(N)$ is generated by $E_{n, n+1}- E_{n+1, n}$, $n=1, \dots, N-1$.

 (ii) The spherical vector in the Gelfand-Zetlin representation is characterized by the set of difference equations
 \beq\label{spheric}
 \left(E_{n, n+1}- E_{n+1, n}\right)\phi_N=0, \; n=1, \dots, N-1.
 \eeq
 (iii) Put
 \beq\label{eq-sph}
 \phi_n(\lambda_{n}, \lambda_{n+1})=\prod_{k=1}^n\prod_{m=1}^{n+1}\Gamma\left(\frac{\lambda_{nk}-\lambda_{n+1,m}}{2i}+\frac{1}{4}\right).
 \eeq
 The function
 \beq
 \phi_N(\lambda)=\prod_{n=1}^{N-1}e^{-\pi(n-1)/2\sum_{j=1}^n \lambda_{nj}}\phi_n(\lambda_{n}, \lambda_{n+1})
 \eeq
 satisfies equations \eqref{spheric}.
 \end{pro}
 While it seems a rather difficult question whether the Gelfand-Zetlin representation may be integrated to a representation of the corresponding Lie group, it is fairly easy to integrate its restriction to the Cartan subalgebra in $\mathfrak{gl}(N)$ which acts by multiplication operators. In this way we get the following explicit formulae for the spherical and Whittaker functions (or, more precisely, of their restrictions to the Cartan subgroup).
\begin{theorem}
The restriction  of the spherical function on $GL(N)$ to the subgroup of diagonal matrices is given (in exponential parametrization) by
\begin{multline}\label{int_phi}
\phi_{\lambda_N}(x_1, \dots, x_N)=\\
\int\limits_{\mathbb{R}^{N(N-1)/2}}e^{i\sum_{n,j=1}^N\left(\lambda_{nj}-\lambda_{n-1,j}\right)x_n}\times\\
\prod_{n=1}^{N-1}\frac{\prod\limits_{k=1}^n\prod\limits_{m=1}^{n+1} \Gamma\left(\frac{\lambda_{nk}-\lambda_{n+1,m}}{2i}+\frac{1}{4}\right)\Gamma\left(\frac{\lambda_{n+1,m}-\lambda_{nk}}{2i}+\frac{1}{4}\right)}{\prod\limits_{s\neq p}\Gamma\left(i\lambda_{np}-i\lambda_{ns}\right)}\prod\limits_{n=1}^{N-1}\prod\limits_{j\leq n} d\lambda_{nj}.
\end{multline}
\end{theorem}
The integrand in \eqref{int_phi} is meromorphic, which allows to deform the integration contour $\mathbb{R}^{N(N-1)/2}\subset\mathbb{C}^{N(N-1)/2}$ and apply the residue theorem. In this way on can easily get the standard asymptotic expansion of the spherical function. A similar formula holds for the Whittaker functions. The formula below gives the Whittaker function
$$
\psi_\lambda(x)=e^{-\rho(x)}\<e^{-\sum_{n=1}^N x_nE_{nn}}w_N, w'_N\>,
$$
which is slightly different from the one discussed in section \ref{rep}; the matrix coefficient
$$
W_\lambda(e^x)=e^{-\rho(x)}\<e^{-\sum_{n=1}^N x_nE_{nn}}w_N, \phi_N\>
$$
used in \eqref{wh} can be computed in a completely similar way and yield once again the Mellin--Barnes  formula for the open Toda wave function.
\begin{theorem}
The restriction  of the Whittaker function on $GL(N)$ to the subgroup of diagonal matrices is given (in exponential parametrization) by
\begin{multline}\label{int_wh}
\psi_{\lambda_N}(x_1, \dots, x_N)=\\
\int\limits_{C}e^{i\sum_{n,j=1}^N\left(\lambda_{nj}-\lambda_{n-1,j}\right)x_n}\times\\
\prod_{n=1}^{N-1}\frac{\prod\limits_{k=1}^n\prod\limits_{m=1}^{n+1} \Gamma\left(\frac{\lambda_{nk}-\lambda_{n+1,m}}{i}\right)}{\prod\limits_{s\neq p}\Gamma\left(i\lambda_{np}-i\lambda_{ns}\right)}\prod\limits_{n=1}^{N-1}\prod\limits_{j\leq n} d\lambda_{nj},
\end{multline}
where the integration contour $C\simeq \mathbb{R}^{N(N-1)/2}\subset \mathbb{C}^{N(N-1)/2}$ is chosen in such a way that
$$
\min\limits_{j}\Im \lambda_{nj}>\max\limits_{m}\Im \lambda_{n+1, m} \q\text{for all}\q n=1,\dots, N-1.
$$
\end{theorem}

Returning back to the definition \eqref{GZ}, one immediately notices its close similarity with the definition of the separation representation  in section \ref{QISM}. Both constructions have intrinsic links with   the deeper aspects of QISM connected with the representation theory of an infinite dimensional algebra, the Yangian $Y(\mathfrak{gl}(N))$ introduced by Drinfeld \cite{Dr}. In the theory of Quantum Groups it is quite common to encounter natural representations of various quantum algebras acting by difference operators.  Remarkably, the definition \eqref{GZ} admits a natural q-deformation which yields a representation of the quantum universal enveloping algebra $U_q(\mathfrak{gl}(N)$  and is the basis for the study of Whittaker and spherical functions in the  q-deformed case\footnote{A q-deformed version of the Mellin--Barnes formula is derived in \cite{KLS}; the  Whittaker functions in the q-deformed case are closely connected with a deep new concept, the modular duality discovered by Faddeev in \cite{Fadd}, and involve a new class of special functions which generalize the Gamma functions.}. However, the emergence of difference operators in the classical setting of semisimple Lie algebras is a totally new phenomenon, which looks both challenging and very promising.
To end up, I would like to make a list of open questions which deserve further study.
\begin{enumerate}
  \item  Integrability problem: is it possible to integrate representations of the universal enveloping algebra by difference operators to the corresponding Lie group. The study of its representations the universal enveloping algebra  by differential operators has a long history which is based on the well-developed spectral theory of differential operators. By contrast, the spectral theory of difference operators of the kind encountered in representation theory (unbounded operators acting in the complex domain) remains almost totally unexplored.

  \item Formulae \eqref{Wh_sol}, \eqref{eq-sph} show that  in the Gelfand--Zetlin representation  both Whittaker and spherical vectors are given by explicit factorized expressions. Is it possible to control in a similar way the intertwining operators for the principal series?
  \item  Both the construction of the principal series representations and the Whittaker theory have natural generalizations to the case of semisimple Lie groups over local fields, such as $\mathbb{Q}_p$, and to the case of global fields. It would be very interesting to construct separation representation and Gelfand--Zetlin representation in these cases.

       \item  Is it possible to generalize separation representation and factorized formulae to other representations of  semisimple Lie groups which are covered by the ``Harish-Chandra philosophy'', i.e. to  general representations induced from the parabolic subgroups? Does the quantum separation of variables technique apply to the study of characters of semisimple Lie groups?
\end{enumerate}
Back in the 1950's and 1960's the fundamental works of Gelfand and Harish-Chandra have created a paradigm which determined the development of the representation theory of  semisimple Lie groups for several decades. The new approach based on the quantum separation of variables announces probably the creation of a new paradigm of no less importance for the future theory.

\end{document}